\documentclass[a4paper,11pt]{article}

\newtheorem{theorem}{Theorem}
\newtheorem{lemma}[theorem]{Lemma}

\newtheorem{corollary}[theorem]{Corollary}
\newtheorem{definition}[theorem]{Definition}
\newtheorem{example}[theorem]{Example}

\newtheorem{proposition}[theorem]{Proposition}
\newtheorem{result}[theorem]{Result}
\newtheorem{observation}[theorem]{Observation}

\newcommand{\bt}{\begin{theorem}}
\newcommand{\et}{\end{theorem}}
\newcommand{\bl}{\begin{lemma}}
\newcommand{\el}{\end{lemma}}
\newcommand{\bp}{\begin{proposition}}
\newcommand{\ep}{\end{proposition}}

\newcommand{\Proof}{ \noindent{\bf Proof:}\quad }
\newcommand{\qed  }{\hfill$\Box$ \smallskip}

\def\qed{\nobreak\hfill\vrule height5pt width4pt
depth2pt\break\medskip}

\def\PG{{\rm PG}}
\def\AG{{\rm AG}}

\usepackage{amssymb}
\usepackage{amsmath}
\usepackage{latexsym}

\usepackage{url}
\usepackage[colorlinks=true,citecolor=black,linkcolor=black,urlcolor=blue]{hyperref}

\newcommand{\doi}[1]{\href{http://dx.doi.org/#1}{\texttt{doi:#1}}}

\title{\bf An extension of the direction problem}
\author{P\'eter Sziklai and Marcella Tak\'ats \thanks{The authors were partially supported by OTKA K 81310, OTKA CNK77780, T\'AMOP, ERC, Bolyai grants.}}

\date{}

\begin{document}

\maketitle

{\small\noindent\MakeUppercase{Notice:} this is the author's version of a work that was accepted for publication in \emph{Discr. Math.} Changes resulting from the publishing process, such as peer review, editing, corrections, structural formatting, and other quality control mechanisms may not be reflected in this document. Changes may have been made to this work since it was submitted for publication. A definitive version was subsequently published in \emph{Discr. Math.}, Volume~\textbf{312}, Issues~12--13 (6~July 2012), Pages 2083--2087. \doi{10.1016/j.disc.2012.02.021}}\\

\begin{abstract}
Let $U$ be a point set in the $n$-dimensional affine space $\AG(n,q)$ over the finite field of $q$ elements and $0\leq k\leq n-2$. In this paper we extend the definition of {\em directions determined by $U$}: a $k$-dimensional subspace $S_k$ at infinity is {\em determined} by $U$ if there is an affine $(k+1)$-dimensional subspace $T_{k+1}$ through $S_k$ such that $U\cap T_{k+1}$ spans $T_{k+1}$. We examine the extremal case $|U|=q^{n-1}$, and classify point sets {\em not} determining every $k$-subspace in certain cases.
\end{abstract}

\section{Introduction}
Let $\AG(n,q)$ and $\PG(n,q)$ denote the $n$-dimensional affine and projective geometries over the finite field of $q$ elements. We will think about $\PG(n,q)$ as $\AG(n,q)\cup H_\infty$, where $H_\infty$ is called the {\em hyperplane at infinity} or the {\em ideal hyperplane}.

Let $U$ be a point set in the affine plane $\AG(2,q)\subset \PG(2,q)$. We say a direction $d$ (i.e. a point at infinity) is {\em determined} by $U$ if there is an affine line with the ideal point $d$ containing at least two points of $U$. It is a well-studied question that how many directions are determined by a point set and how the ``interesting'' point sets $U$ look like, see \cite{Balldirec}, \cite{BBBSS}. It seems to be natural to extend this definition in the following way:\\

\begin{definition}
Let $U\subset \AG(n,q)\subset \PG(n,q)$, and $k$ be a fixed integer, $k \leq n-2$. We say a subspace $S_k$ of dimension $k$ in $H_{\infty}$ is {\em determined} by $U$ if there is an affine subspace $T_{k+1}$ of dimension $k+1$, having $S_k$ as its hyperplane at infinity, containing at least $k+2$ affinely independent points of $U$ (i.e. spanning $T_{k+1}$).
\end{definition}

The question is that for a fixed $k$ how many subspaces of dimension $k$ are determined by $U$, or how large the point set can be, if it does not determine all of the $k$-subspaces at infinity, and can we say something about the structure of $U$ in that case. Consider a subspace $S_k$ of dimension $k$ at infinity. There are $q^{n-k-1}$ (pairwise disjoint) affine subspaces of dimension $k+1$ with the ideal hyperplane $S_k$. Suppose that $S_k$ is not determined, and consider an affine $(k+1)$-subspace $T_{k+1}$ through $S_k$. Then the points of the intersection of $U$ and $T_{k+1}$ are contained in some subspace of dimension $k$ since $S_k$ is undetermined and so $T_{k+1}$ is not spanned by them. An affine subspace of dimension $k$ consists of $q^k$ points, so (any) $T_{k+1}$ can contain at most $q^k$ points of $U$. This implies $|U| \leq q^{n-1}$ if it does not determine all the $k$-subspaces at infinity.\\

\section{The extremal case}

We will examine point sets of the extremal cardinality, so let $|U|=q^{n-1}$. Our aim is to find out a (strong) structure of $U$, if there are relatively many undetermined subspaces at infinity.\\
Let $|U|=q^{n-1}$. Suppose that there exists an undetermined $k$-subspace $S_k \subset H_{\infty}$. Then in each of the $q^{n-k-1}$ affine $(k+1)$-subspaces whose ideal hyperplane is $S_k$ there lie exactly $q^k$ points of $U$ constituting one complete $k$-subspace.\\
We can construct point sets not determining many subspaces in the following way:\\
Let $U_0\subset \AG(m,q)\subset \PG(m,q)$, $|U_0|=q^{m-1}$. Denote by $N_l^0$ the set of the $l$-subspaces in $H_{\infty}$ which are not determined by $U_0$. We can embed $U_0$ into $\AG(n,q)$, where $n>m$. Consider a subspace $V$ on the ideal hyperplane of $\AG(n,q)$, $\dim V=v$, where $v=n-m-1$, completely disjoint from the original $m$-dimensional space. We construct a cone (cylinder) with base $U_0$ and vertex $V$ such that we take the union of the affine subspaces spanned by a point of $U_0$ and $V$. These subspaces are of dimension $n-m$, and in that way we get a point set $U$ in $\AG(n,q)$, $|U|=q^{n-1}$.\\
We will show that there are ``many'' subspaces at infinity which are not determined by $U$, and we characterize them.

Denote by $N_r$ the set of subspaces of dimension $r$ on the ideal hyperplane which are not determined by $U$.

\begin{proposition}
\label{cone}
Let $U_0\subset \AG(m,q)\subset \PG(m,q)$, $|U_0|=q^{m-1}$ embedded into $\AG(n,q)$, let $V \subseteq  H_\infty$ be a subspace in $PG(n,q)$ completely disjoint from $\PG(m,q)$, $\dim V=n-m-1$, let $U\subset \AG(n,q)\subset \PG(n,q)$ be the cone (cylinder) constructed as above with base $U_0$ and vertex $V$. Let $W\subseteq  H_\infty$ in $\PG(n,q)$, $\dim W=r$. Then $W$ is non-determined, i.e. $W\in N_{r}$ if and only if, after projecting  $W$ from $V$ to the $m$-space, the projected image $W_0=\PG(m,q)\cap \langle V,W\rangle$ is non-determined by $U_0$.

\end{proposition}

\Proof
For the proof let $s=\dim(W\cap V)$ and $r_0=\dim W_0=r-s-1$. (Note that $W$ and
$W_0$ are contained in the hyperplane at infinity.) By projection we always
mean projection from the center $V$ (see below).

Suppose that $W$ is non-determined. Let $L_0\supset W_0$ be any of the affine
$r_0+1$-dimensional subspaces in $\PG(m,q)$ through $W_0$, we have to prove
that $L_0\cap U_0$ is a complete affine $r_0$-space.
Take any affine point $P$ from the subspace $\langle V, L_0\rangle$, then
define $L=\langle P, W\rangle$. Now $L_0$ is the projected image of $L$,
i.e. $L_0=\PG(m,q)\cap\langle V,L\rangle$. As $\dim L = r+1$ and $W$ is
non-determined, $A=L\cap U$ must be a complete affine $r$-space. Hence $L_0\cap
U_0$ is the projected image of it, i.e. a complete affine $r_0$-space.
(We have used the fact that $\dim(\bar{A}\cap V)=s=\dim(W\cap V)$, where
$\bar{A}$ is the projective closure of $A$. Indeed, if $\dim(\bar{A}\cap V)$
were larger then $\dim(L\cap U)$ would be larger accordingly.)
So $W_0$ is non-determined.

On the other hand, suppose that $W_0$ is non-determined. Then for any
$r+1$-dimensional affine subspace $L$ through $W$, and for its projected image
$L_0$, the intersection $L_0\cap U_0$ is identical to the projected image of
$L\cap U$. As  $L_0\cap U_0$ is a complete $r_0+1$-dimensional affine subspace
through $W_0$, we have that $L\cap U$ must be $\langle V, L_0\rangle \cap L$,
i.e. a complete  $r$-dimensional affine subspace through $W$.
\qed

For a given affine point set we can examine determined subspaces in $H_\infty$ of different dimensions. We are trying to find out a hierarchy of the determined subspaces of different dimensions. In \cite{StSz} the following theorem was proved in the classical case:

\begin{result}
  \label{StormeSziklai}
Let $U\subset \AG(n,q)\subset \PG(n,q)$, $|U|=q^{n-1}$ and let $D \subseteq H_{\infty}$ be the set of directions determined by $U$. Then $D$ is the union of some complete $(n-2)$-dimensional subspaces of $H_{\infty}$.
\qed
\end{result}

We find an analogous situation when determining higher dimensional subspaces:

\begin{observation}
Let $U\subset \AG(n,q)\subset \PG(n,q)$, $|U|=q^{n-1}$ and $k$ be a fixed integer, $k \leq n-3$. If there is a subspace $S_k$ of dimension $k$, $S_k \subset H_{\infty}$ determined by $U$, then there is a subspace $S_{k+1}$ of dimension $k+1$ in $H_{\infty}$, $S_k \subset S_{k+1}$, which is determined by $U$.
\end{observation}

\Proof
Since $S_k $ is determined by $U$, there exists an affine $(k+1)$-dimensional subspace $A_{k+1}$ with the ideal hyperplane $S_k$ which contains at least $k+2$ linearly independent points from $U$. There will be at least one $(k+2)$-dimensional affine subspace $A_{k+2}$, $A_{k+1} \subset A_{k+2}$ containing at least one more point from $U$. So $A_{k+2}$ contains at least $k+3$ linearly independent points from $U$, so its ideal hyperplane will be determined.
\qed

\begin{corollary}
Analogously to {\bf Result \ref{StormeSziklai}}, through any determined $k$-subspace $S_k$, there exists a determined $(n-2)$-subspace.
\qed
\end{corollary}

\begin{proposition}
Let $U\subset \AG(n,q)\subset \PG(n,q)$, $|U|=q^{n-1}$ and $k$ be a fixed integer, $k \leq n-3$. If there is a subspace $V_{n-2}$ of dimension $n-2$, $V_{n-2} \subset H_{\infty}$ such that all of the $k$-dimensional subspaces of $V_{n-2}$ are determined by $U$ then $V_{n-2}$ is determined by $U$ as well.
\end{proposition}

\Proof
Suppose that there is a subspace $V_{n-2}$ of dimension $n-2$, $V_{n-2} \subset H_{\infty}$ which is not determined by $U$. We will show that there is a subspace $V_k$, $\dim V_k=k$, $V_k\subset V_{n-2}$ which is not determined by $U$.
Each of the $q$ affine $(n-1)$-subspaces whose ideal hyperplane is $V_{n-2}$ contains precisely a whole $(n-2)$-subspace from $U$. Such an $(n-2)$-subspace has an ideal hyperplane of dimension $n-3$ contained in $V_{n-2}$, so these $(n-3)$-subspaces cannot cover all the points of $V_{n-2}$. Consider an uncovered point $P$, and take a subspace $V_k\subset V_{n-2}$, $\dim V_k=k$ containing $P$. $V_{n-2}$ and each of the affine subspaces of dimension $k+1$ with the ideal hyperplane $V_k$ span a subspace of dimension $n-1$. Such an $(n-1)$-subspace contains precisely a whole $(n-2)$-subspace from $U$. The intersection of the $(n-2)$-subspace from $U$ and the $(k+1)$-subspace is a subspace of dimension $k$ (since this $(n-2)$-subspace cannot contain the $(k+1)$-subspace because $P \in V_k$.) So each of the $(k+1)$-spaces with the ideal hyperplane $V_k$ contains precisely a whole $k$-space from $U$ which means that $V_k$ is not determined by $U$.
\qed

\section{The 3-dimensional case}

The question is whether we can give other constructions for $U$ with many undetermined subspaces at infinity, or can we say something about the point set if there are less undetermined subspaces. In certain dimensions we found that in case of a few undetermined subspaces the point set has a strong structure.\\

Let $U\subset \AG(3,q)\subset \PG(3,q)$. We say a line $\ell \subset H_\infty$ is determined by $U$ if there is an affine plane with the ideal line $\ell$ containing at least three points of $U$ which are not collinear. Let $|U|=q^2$. Suppose that there exists an undetermined line $\ell$. Then each of the $q$ affine planes whose ideal line is $\ell$ contains precisely a complete line of $U$.

\begin{theorem}
\label{lindet}
Let $U\subset \AG(3,q)\subset \PG(3,q)$, $|U|=q^2$. Let $L$ be the set of lines in $H_\infty$ determined by $U$ and put $N$ the set of non-determined lines. Then one of the following holds:\\

\noindent{\bf a)} $|N| = 0$, i.e.  $U$ determines all the lines of $H_\infty$;\\
\noindent{\bf b)} $|N| = 1$ and then there is a parallel class of affine planes such that $U$ contains one (arbitrary) complete line in each of its planes;\\
\noindent{\bf c)} $|N| = 2$ and then $U$ together with the two undetermined lines in $H_\infty$ form a hyperbolic quadric or $U$ contains $q$ parallel lines ($U$ is a cylinder);\\
\noindent{\bf d)} $|N| \geq 3$ and then $U$ contains $q$ parallel lines ($U$ is a cylinder).\\

\end{theorem}

\Proof
If all the points of $U$ are contained in a plane, then $U$ determines only one line, the line in $H_\infty$ of the plane containing $U$, so $|L| = 1$. It is a special case of {\bf d)}.\\
So from now on suppose that $U$ is not contained in a plane.\\
Denote by $f_\infty$ the ideal point of an affine line $f$, and let $\bar f=f \cup f_\infty$.\\
Throughout this proof let $|N|\geq 2$, $\ell_1, \ell_2 \subset H_\infty$ be undetermined lines, $M$ the intersection point of $\ell_1$ and $\ell_2$.

\begin{lemma}
\label{1meet2}
Let $f$ and $g$ be affine lines contained in $U$, $\bar f$ intersecting $\ell_1$, $\bar g$ intersecting $\ell_2$. Then $f$ and $g$ meet each other or their ideal point is $M$.
\end{lemma}

\Proof
Suppose to the contrary that $\bar f$ and $\bar g$ are skew lines, and at least for one of them, e. g. $\bar f$, $M \notin \bar f$. Then $\bar f$ intersects the plane $s$ determined by $\bar g$ and $\ell_2$ at an affine point $P$ not contained in $g$. Then $P$ and $g$ would determine $\ell_2$, contradiction.
\qed

\begin{lemma}
\label{l1lines}
If there is an affine line $f$ contained in $U$, $\bar f$ intersecting $\ell_1$ at $M$, then the ideal point of each of the $q$ affine lines contained in $U$ whose ideal points are in $\ell_1$ is $M$.
\end{lemma}

\Proof
Suppose that there exists a line $g$ contained in $U$, $\bar g$ intersecting $\ell_1$ not at $M$. Since $\bar f$ intersects $\ell_1$ at $M$, it intersects $\ell_2$ as well, so by {\bf Lemma \ref{1meet2}} $f$ meets $g$ at an affine point $P$. Then $f$ and $g$ would determine $\ell_1$, contradiction.
\qed

\begin{lemma}
\label{2parlin}
If there exist $f$ and $g$ parallel affine lines contained in $U$, $\bar f, \bar g$ intersecting $\ell_1$, then their ideal point is $M$.
\end{lemma}

\Proof
Suppose to the contrary that $\bar f$ and $\bar g$ meet $\ell_1$ not at $M$. Consider a point $P \in U$ not contained in the plane $\langle f,g\rangle$, and the plane $s$ determined by $P$ and $\ell_2$. (If $\nexists P \in U$, $P\notin\langle f,g\rangle$ then all the points of $U$ would be contained in a plane.) $s$ cannot be parallel to $\bar f$ and $\bar g$ as then they would intersect in the ideal hyperplane (i.e. in $\ell_2$). So $s$ intersects $f$ and $g$ at two different affine points. The two intersection points and $P$ cannot be collinear as $P\notin\langle f,g\rangle$, so they span $s$ determining $\ell_2$, contradiction.
\qed

\begin{corollary}
\label{allpar}
If there exist $f$ and $g$ parallel affine lines contained in $U$, $\bar f, \bar g$ intersecting $\ell_1$, then $U$ consists of $q$ parallel lines whose ideal point is $M$.
\end{corollary}

\Proof
Since $f$ and $g$ are parallel, $\bar f$ and $\bar g$ meet at $M$ due to {\bf Lemma \ref{2parlin}}. Then, by {\bf Lemma \ref{l1lines}}, for all the lines $h$ contained in $U$ whose ideal points are in $\ell_1$, $\bar h$ will intersect $\ell_1$ at $M$, which means that they are all parallel.
\qed

\begin{corollary}
\label{undetint}
If $|N|> 2$, $\ell_1 \subset H_\infty$ is an undetermined line, and there exist $f$, $g$ parallel lines contained in $U$, $\bar f$, $\bar g$ intersecting $\ell_1$, then all the undetermined lines in $H_{\infty}$ intersect at the same point (and it is the ideal point of the lines of $U$).
\end{corollary}

\Proof
Let $\ell_2$, $\ell_3$ be undetermined lines, $M$ the intersection point of $\ell_1$ and $\ell_2$, $K$ the intersection point of $\ell_1$ and $\ell_3$. By { \bf Lemma \ref{2parlin}}, $\bar f$ and $\bar g$ meet $\ell_1$ at $M$ and also at $K$.
\qed

\begin{corollary}
\label{allskew}
If there exist $2$ lines contained in $U$ whose ideal points are in $\ell_1$ which are not parallel, then all the $q$ lines whose ideal points are in $\ell_1$ have $q$ different ideal points (and none of them is $M$).
\end{corollary}

\Proof
By {\bf Corollary \ref{allpar}}, if there exist two parallel lines, then all of them are parallel. So if they are not all parallel, then all of them are pairwise skew. It means that their ideal points in $\ell_1$ are pairwise different. So there is a line $f \subset U$ such that $\bar f$ intersects $\ell_1$ not at $M$. Then by {\bf Lemma \ref{l1lines}} there cannot exist a line $g \subset U$ such that $\bar g$ intersects $\ell_1$ at $M$.
\qed

So we have two different cases: {\bf (1)} all the $q$ lines meeting $\ell_1 \in N$ are parallel and intersect the ideal hyperplane at the intersection point of the undetermined lines or {\bf (2)} all of them are pairwise skew.\\
In the first case $U$ consists of $q$ (say) vertical lines (forming a cylinder), and the undetermined lines intersect each other at the ideal point of the vertical lines. Then every ``horizontal" (affine) plane (intersecting the vertical lines) contains $q$ points of $U$. Consider the directions not determined by these $q$ points on a ``horizontal" plane.  The undetermined lines are exactly the lines connecting the ideal point of the vertical lines and an undetermined direction.\\
It is exactly the construction we saw in {\bf Proposition \ref{cone}}: the base of the cone is a point set of cardinality $q$ contained in a horizontal plane, and the vertex of the cone is $M$.\\
There is a special case of this first case: the $q$ vertical lines can be co-planar. Then $U$ is contained in a plane. Then $U$ determines only one line, so $|L| = 1$.\\

In the second case let $\ell_1$, $\ell_2$ be undetermined lines. By {\bf Corollary \ref{allskew}} we know that the points of $U$ form $q$ lines whose ideal points are in $\ell_1$, and these ideal points are pairwise different, so the lines are skew.
It also holds for $\ell_2$, so $U$ forms $q$ skew lines whose ideal points are in $\ell_2$, and these ideal points are pairwise different. Consider $3$ of the skew lines whose ideal points are in $\ell_1$ and denote them by $f$, $g$ and $h$.
By {\bf Lemma \ref{1meet2}} all the $q$ lines whose ideal points are in $\ell_2$ intersect $f$, $g$ and $h$ at one-one point. The intersection points are different as the $q$ lines are pairwise skew. So the $q$ lines whose ideal points are in $\ell_2$ form a $q$-regulus of $f$, $g$ and $h$, and with $\ell_1$ it is a $(q+1)$-regulus of $\bar f$, $\bar g$ and $\bar h$.
The same holds for the lines $f'$, $g'$ and $h'$ whose ideal points are in $\ell_2$: the $q$ lines whose ideal points are in $\ell_1$ together with $\ell_2$ form a $(q+1)$-regulus of $\Bar {f'}$, $\Bar {g'}$ and $\Bar {h'}$. It means that the $q^2$ points of $U$ and the undetermined lines $\ell_1$ and $\ell_2$ form a hyperbolic quadric and {\bf Theorem \ref{lindet}} is proved.\qed

\section{More quadrics}

So we have the complete characterization in $3$ dimensions. The question is whether we can have similar results in higher dimensions. In $\PG (n,q)$ we will show that if we try to find a point set $U$, $|U|=q^{n-1}$ which does not determine all the subspaces of a certain dimension then the former examples will occur:\\
If the point set $U$ corresponds to the affine part of a nonsingular quadric for which the ideal hyperplane is a tangent hyperplane, or $U$ is a cone constructed as in {\bf Proposition \ref{cone}}, then there will be some undetermined subspaces. The open question is whether there are other point sets not determining more than one subspaces or nonsingular quadrics and cones are the only examples.\\

Let $U\subset \AG(n,q)\subset \PG(n,q)$, $|U|=q^{n-1}$ be the affine part of a nonsingular quadric for which the ideal hyperplane is a tangent hyperplane (i.e. the intersection of the quadric and the ideal hyperplane is a cone based on an $(n-2)$-dimensional quadric of the same character). Denote by $g$ the projective index of the quadric, i.e. the dimension of the generators, the subspaces of maximum dimension contained in the quadric. Denote by $w$ the character of the quadric, $w=1,\ 2,\ 0$ for a parabolic, hyperbolic, elliptic quadric, respectively.

\begin{proposition}
If $U\subset \AG(n,q)\subset \PG(n,q)$, $|U|=q^{n-1}$ is defined as above then the undetermined $g$-dimensional subspaces in $H_\infty$ are exactly the generators contained in the intersection of the quadric and the ideal hyperplane, except the cases when $n=2$ and $q$ is even, or $n=4$ and $q=2$, or $n=5$ and $q=2$ and the quadric is elliptic.
\end{proposition}

\Proof
There are $q^{n-1}$ affine points in the quadric. Denote by $G$ a generator of the quadric contained in $H_\infty$. There are $q^{n-g-1}$ affine $(g+1)$-dimensional subspaces with the ideal hyperplane $G$. Each of these $(g+1)$-subspaces contains at most $q^g$ points from the quadric, otherwise the whole $(g+1)$-subspace would be contained in the quadric, which is a contradiction. Since there are $q^{n-1}$ affine points in the quadric, each $(g+1)$-subspace contains exactly $q^g$ points. The intersection of a $(g+1)$-subspace and the quadric is a quadric of dimension $g+1$, containing $G$, a subspace of dimension $g$, which is a linear factor. So the rest of the intersection has to be a $g$-subspace as well. It means that each of the $(g+1)$-dimensional affine subspaces with the ideal hyperplane $G$ contains exactly an affine $g$-dimensional subspace from $U$, so $G$ is undetermined.\\

On the other hand consider a $g$-dimensional subspace $G'$ in $H_\infty$ not contained in the quadric. Suppose it is undetermined. It means that each of the $q^{n-g-1}$ \ \ $(g+1)$-dimensional affine subspaces with the ideal hyperplane $G'$ contains exactly a $g$-dimensional subspace from $U$ (i.e. a generator). If the quadric is hyperbolic, such a $(g+1)$-dimensional subspace should contain one more generator (from the other system), contradiction. If the quadric is elliptic or parabolic, the ideal hyperplanes of these generators are in $G'$ and they are of dimension $g-1$. Since $G'$ is not a generator, it can contain at most $2$ \ \ $(g-1)$-subspaces from the quadric. This implies that there is a $(g-1)$-subspace in the intersection of $G'$ and the quadric, which is the ideal hyperplane of at least $\frac{q^{n-g-1}}{2}$ generators. The number of generators through a fixed $(g-1)$-subspace is the following:

$\varrho (g-1,n,w)= (q^{2-w}+1)\cdot(q^{2-w+1}+1)\cdot ... \cdot(q^{\frac{n-2g+1-w}{2}}+1)$.

For an elliptic quadric $\varrho (\frac{n-5}{2},n,0)= q^2+1 < \frac{q^{\frac{n+1}{2}}}{2}$, if $n \geq 5$ and $q \geq 3$\ or $n \geq 7$;\\
for a parabolic quadric $\varrho (\frac{n-4}{2},n,1)= q+1 < \frac{q^{\frac{n}{2}}}{2}$, if $n \geq 4$ and $q \geq 3$\ or $n \geq 6$, contradiction.

If $n=3$, for an elliptic quadric it is clear that through a point in $H_\infty$ not contained in the quadric there are $q^2$ affine lines, and $q+1$ of them are tangents to the quadric, all the others are secants or skew lines. Since there are $q^2$ affine points in the quadric, the point at infinity will be determined.

Similarly, if $n=2$, $q$ odd, through a point in $l_\infty$ not contained in the quadric, there is only one affine line which is tangent to the quadric, all the others are secants or skew lines, so the point will be determined.

There are exceptional cases. If $n=2$ and $q$ is even, the nucleus of the parabolic quadric is a not determined point in the ideal line, but it is not contained in the quadric. For a parabolic quadric in $\PG(4,2)$ and for an elliptic quadric in $\PG(5,2)$ there exist undetermined lines in $H_\infty$ which are not generators.
\qed


\begin{thebibliography}{99}

\bibitem{Balldirec} {\sc S.~Ball}, The number of directions determined by a function
over a finite field,
{\em J. Combin. Th. Ser. A}, {\bf 104} (2003), 341--350.


\bibitem{BBBSS} {\sc A. Blokhuis, S. Ball, A. Brouwer, L. Storme and
T. Sz\H onyi}, On the number of slopes determined by a function on
a finite field, {\em J. Comb. Theory Ser. (A)} {\bf 86} (1999),
187--196.


\bibitem{StSz}
 {\sc L. Storme and P. Sziklai}, Linear point sets and R\'edei type $k$-blocking sets in $\PG(n,q)$,
  {\em J. Alg. Comb.} {\bf 14} (2001), 221-228.

\bibitem{HirT}
{\sc J. W. P. Hirschfeld and J. A. Thas}, {\em General Galois geometries}, Clarendon Press, Oxford (1991)

\end{thebibliography}
\end{document}